\documentclass[letterpaper,12pt]{article}
\title{The Maximum Likelihood Degree}

\author{Fabrizio Catanese, Serkan Ho\c{s}ten, \\
Amit Khetan, and Bernd Sturmfels}

\date{}

\usepackage{latexsym,array,delarray,amsthm,amssymb,epsfig}
\usepackage[sumlimits]{amsmath}

\theoremstyle{plain}
\newtheorem{thm}{Theorem}
\newtheorem{lemma}[thm]{Lemma}
\newtheorem{prop}[thm]{Proposition}
\newtheorem{cor}[thm]{Corollary}

\theoremstyle{definition}
\newtheorem{df}[thm]{Definition}
\newtheorem{ex}[thm]{Example}
\newtheorem{rem}[thm]{Remark}

\newcommand{\ra}{\ensuremath{\,\rightarrow}\,}

\def\eea{\end{eqnarray*}}
\def\bea{\begin{eqnarray*}}

\def\X{{\mathcal{X}}}

\newcommand{\C}{\ensuremath{\mathbb{C}}}

\newcommand{\hol}{\ensuremath{\mathcal{O}}}

\newcommand{\PP}{\ensuremath{\mathbb{P}}}

\newcommand{\FF}{\ensuremath{\mathcal{F}}}
\newcommand{\EE}{\ensuremath{\mathcal{E}}}
\newcommand{\HHH}{\ensuremath{\mathcal{H}}}

\newcommand{\zz}{\mathbb{Z}}
\newcommand{\nn}{\mathbb{N}}
\newcommand{\pp}{\mathbb{P}}

\newcommand{\rr}{\mathbb{R}}
\newcommand{\cc}{\mathbb{C}}

\begin{document}
\maketitle
\begin{abstract}

Maximum likelihood estimation in statistics leads to the
problem of maximizing a product of powers of polynomials.
We study the algebraic degree of the critical equations
of this optimization problem. This degree is related to
the number of bounded regions in the  corresponding
arrangement of hypersurfaces, and to the Euler
characteristic of the complexified complement.
Under suitable hypotheses,  the maximum likelihood degree
 equals the top Chern class of
a sheaf of logarithmic differential forms.
Exact formulae
in terms of degrees and Newton polytopes
are given for polynomials with generic coefficients.
     \end{abstract}

\section{Introduction}

In algebraic statistics \cite{GSS,PS,PRW}, a {\em model for discrete data} is a
    map ${\bf f}  : \rr^d  \rightarrow \rr^n$
whose coordinates $f_1,\ldots,f_n$ are polynomial functions
in the parameters $ (\theta_1,\ldots,\theta_d)=:\theta  $.
The parameter vector $\theta$ ranges over an open  subset
$\mathcal{U}$ of $\rr^d$ such that $f(\theta)$
lies in the positive orthant $\rr^n_{> 0}$.
The image  $f(\mathcal{U})$ represents a
family of probability distributions on
an $n$-element state space, provided
we make the extra assumption that
$\,f_1 + \cdots + f_n -1\,$
is the zero polynomial.

A given data set is a vector
$\,u = (u_1,\ldots,u_n)\,$ of positive integers.
The problem of {\em maximum likelihood estimation}
is to find parameters $\theta$ which best explain the data $u$.
This leads to the following optimization problem:
\begin{equation}
\label{mle}
    {\rm Maximize} \,\,\,
f_1(\theta)^{u_1} f_2(\theta)^{u_2} \cdots f_n(\theta)^{u_n} \quad
\hbox{subject to $\,\,\,\theta \in \mathcal{U}$} .
\end{equation}
Under suitable assumptions we have an optimal solution
$\hat{\theta}$ to the problem (\ref{mle}), which is an algebraic
function of the data $u$. Our goal is to compute the degree of
that algebraic function. We call this number the {\em maximum
likelihood degree} of the model ${\bf f}$. Equivalently, the {\em
ML degree} is the number of complex solutions to  the critical
equations of (\ref{mle}), for a general data vector $u$. In this paper
we
   prove results of the following form:

\begin{thm} \label{thm1}
Let $f_1,\ldots,f_n$ be polynomials
of degrees $b_1,\ldots,b_n$ in $d$ unknowns. If
   the maximum likelihood degree of the model
   $\,{\bf f} = (f_1,\ldots,f_n)$ is finite then it
is less than or equal to the coefficient of $z^d$ in the generating function
\begin{equation}
\label{thm1formula}
    \frac{ (1-z)^d }{ (1-z b_1)  (1-z b_2)  \cdots (1-z b_n) }.
\end{equation}
Equality holds if the coefficients of the
polynomials $f_i$
are sufficiently generic.
\end{thm}

As an example, consider a model given by $n=4$
quadratic polynomials in $d=2$ parameters.
The solution to (\ref{mle}) satisfies the two critical equations
$$
\frac{u_1}{f_1} \frac{\partial f_1}{\partial \theta_1}
+\frac{u_2}{f_2} \frac{\partial f_2}{\partial \theta_1}
+\frac{u_3}{f_3} \frac{\partial f_3}{\partial \theta_1}
+\frac{u_4}{f_4} \frac{\partial f_4}{\partial \theta_1}
\,=\, \frac{u_1}{f_1} \frac{\partial f_1}{\partial \theta_2}
+\frac{u_2}{f_2} \frac{\partial f_2}{\partial \theta_2}
+\frac{u_3}{f_3} \frac{\partial f_3}{\partial \theta_2}
+\frac{u_4}{f_4} \frac{\partial f_4}{\partial \theta_2}
\,=\, 0 .
$$
If the $f_i$'s are general quadrics then
these equations have
$25$ complex solutions. The formula for
the maximum likelihood degree in Theorem \ref{thm1} gives
$$
    \frac{(1-z)^2}{(1-2z)^4} \,= \, 1 + 6 z + \underline{25} z^2 +
88 z^3 +  280 z^4 + \cdots.
$$
For special quadrics $f_i$, the ML degree
can be much lower than $25$.
A familiar example is the {\em independence model for
two binary random variables}:
\begin{equation}
\label{2by2table}
    f_1 \,=\, \theta_1 \theta_2, \,\,\,
    f_2 \,=\, (1- \theta_1) \theta_2, \,\,\,
    f_3 \,=\, \theta_1 (1-\theta_2), \,\,\,
    f_4 \,=\, (1-\theta_1) (1-\theta_2).
\end{equation}
Here the ML degree is only one because
the maximum likelihood estimate $\hat{\theta}$ is a rational
function (= algebraic function of degree one)
of the data $u$:
$$
\hat{\theta}_1 = \frac{u_1+u_3}{u_1+u_2+u_3+u_4} \quad \hbox{and} \quad
\hat{\theta}_2 = \frac{u_1+u_2}{u_1+u_2+u_3+u_4}.
$$

This paper is organized as follows.
In Section 2 we present the algebraic geometry
for studying critical points of a
rational function $\,f = f_1^{u_1} \cdots f_n^{u_n}\,$
on an irreducible projective variety $X$.
The critical equations $\, d log(f) = 0 \,$
are interpreted as sections of the sheaf
$\,\Omega^1( log \, D )\,$
of $1$-forms with logarithmic singularities
along the divisor $D$ defined by $f$.
In Theorem  \ref{chernclassthm},
   we show that if $D$ is a global normal crossing
divisor then the ML degree equals the degree of the top Chern class of
$\,\Omega^1(log \, D) $. If  $X$ is projective $d$-space then this
leads to Theorem \ref{mle}. In Section 3 we study the case when
$X$ is a smooth toric variety, and we derive a formula for the ML
degree when the $f_i$'s are Laurent polynomials which are generic
relative to their Newton polytopes. For instance, Example
\ref{rectangles} shows that the ML degree is $13$  if we replace
(\ref{2by2table}) by 
$$\,f_i \quad = \quad \alpha_i + \beta_i \theta_1 +
\gamma_i \theta_2 + \delta_i \theta_1 \theta_2
\quad \qquad (i=1,2,3,4). $$

Section 4 is concerned with the relationship
of the ML degree to the bounded regions of the complement
of $\{f_i = 0 \}$ in $\rr^d$. The number of these
regions is a lower bound to the number of real solutions of
the critical equations, and therefore a lower bound to
the ML degree. We show that for plane quadrics all
three numbers can be equal. However, for other
combinations of plane curves the ML degree and the
number of bounded regions diverge, and we prove
a tight upper bound on the latter in Theorem  \ref{thm:viro}.
Also, following work
of Terao \cite{T} and Varchenko \cite{V},
we show in Theorem \ref{thm-linear} that the ML degree
coincides with the number of bounded regions of the
arrangement of hyperplanes $\{f_i = 0\}$ when the $f_i$'s are
(not necessarily generic) linear forms.

Section 5 revisits the ML degree for toric varieties, replacing the
smoothness assumption by a much milder condition. Theorem
\ref{toricformula4} gives a purely combinatorial formula
for the ML degree in terms of the Newton polytopes of the
polynomials $f_i$. This section also discusses how resolution of
singularities can be used to compute the ML degree for nongeneric
polynomials.

Section 6 deals with topological methods for determining the ML degree.
Theorem \ref{euler} shows that, under certain restrictive hypotheses, 
it coincides with the Euler characteristic 
of the complex manifold $X \backslash D$,  and
Theorem \ref{semiconti} offers a general version of the 
semi-continuity principle which underlies the inequality
in Theorem \ref{thm1}. In Section 7 we relate the ML degree to
the sheaf of logarithmic vector fields along $D$,
which is the sheaf dual to $\Omega^1(log D)$. 

This paper was motivated by recent appearances  of the concept of
ML degree in statistics and computational
   biology. Chor, Khetan and Snir \cite{CKS}
   showed that the ML degree of a phylogenetic model equals $9$,
   and Geiger, Meek and Sturmfels \cite{GMS} proved that an
   undirected graphical model
   has ML degree one if and only if it is
decomposable. The notion of ML degree also makes sense for
certain parametrized models for continuous data:
   Drton and Richardson \cite{DR} showed that the ML degree
   of a Gaussian graphical model equals $5$,
and Bout and Richards \cite{BR} studied the ML degree
of certain mixture models.
   The ML degree always provides an upper bound
on the number of local maxima of the likelihood function.
  Our ultimate hope is that
a better understanding of the ML degree will
lead to the development of custom-tailored algorithms
  for solving the critical equations $\,dlog(f) = 0$.  There is a
need for such new algorithms, given that methods currently used in
statistics (notably the {\em EM-algorithm}) often produce only local
maxima in~(\ref{mle}).

\section{Critical Points of Rational Functions}

In this section we work in the following
general set-up of algebraic geometry.
Let $X$ be a complete factorial algebraic variety
over the complex numbers $\C$.
We also assume that $X$ is irreducible of dimension $d \geq 1$.
In applications to statistics, the variety $X$ will
often be a  smooth projective toric variety.

Suppose that $f \in \C (X)$  is a rational function on $X$.
Since $X$ is factorial, the local rings $\hol_{X,x}$ are
unique factorization domains. This means that
the function $f$ has a global factorization
which is unique up to constants:
\begin{equation}
\label{globfac}
   f \quad = \quad F_1^{u_1} F_2^{u_2} \cdots F_r^{u_r}.
\end{equation}
Here $F_i$ is a prime section of an invertible sheaf $\hol_X(D_i)$
where $ D_i $ is the divisor on $X$ defined by $F_i$.
In our applications we usually assume that 
$r \geq n$  where $n$ is the number 
considered in the Introduction.
For instance, if $f_1,\ldots,f_n$ are polynomials
and  $X = \PP^d$ then $r = n+1$; namely,
  $F_1, \ldots, F_n$ are the homogenizations of $f_1,
\ldots, f_n$ using $\theta_0$, and $F_{n+1} = \theta_0$ (see the proof of
Theorem \ref{thm1} for details).

By (\ref{globfac}), we can write
the divisor of the rational function $f$ uniquely as
$$ div(f) \quad = \quad \sum_{i=1}^{r} u_i D_i, $$
where the $u_i$'s are (possibly negative) integers.
Let $D$ be the reduced  union of the
codimension one subvarieties $D_i \subset X$, or, as a divisor,
$\,D : = \Sigma_{i=1}^{r}  D_i$.

We are interested in computing the critical points
of the rational function $f$ on the open set
$\, V := X \backslash D\,$ complementary to the divisor $D$.
Especially,
we wish to know the number of critical points,
counted with  multiplicities.

A critical point is by definition a point $x \in X$ where
the differential $1$-form $df$ vanishes.
If $x$ is a smooth point on $X$, and
$x_1, \dots, x_d$  are local coordinates,
then  $\, df = \Sigma_{j=1}^{d}
(\partial f/ \partial x_j ) dx_j$.
Hence $x$ is a critical point of $f$ if and only if
\begin{equation}
\label{criticalequations}
   \frac{\partial f}{\partial  x_1}  \,\, = \,\,
\frac{\partial f}{\partial  x_2}  \,\, = \,\,
\cdots  \,\, = \,\,
\frac{\partial f}{\partial  x_d}   \,\, = \,\, 0 .
\end{equation}
We next rewrite the critical equations
(\ref{criticalequations}) using
the factorization (\ref{globfac}).
Around each point $x \in X$, we may choose a local
trivialization for the sheaf $\hol_X(D_i)$
and  express $F_i$ locally by a regular function.
By slight abuse of notation, we denote that regular
function also by $F_i$. For instance, if $X = \PP^d$
then this means replacing the homogeneous
polynomial $F_i$ by a dehomogenization.

Since $f$ has neither  zeros nor poles
on the open set $V$, the vanishing of $ df$ is equivalent
to the vanishing of  the logarithmic derivative
\begin{equation}
\label{logderive}
   d log (f) \,\,\, = \,\,\, \frac{df}{f}
\quad = \quad
\sum_{i=1}^{r} u_i  \ d log (F_i)
\quad = \quad
\sum_{i=1}^{r} u_i  \ \frac{dF_i}{F_i}.
\end{equation}

We now recall some classical definitions and results
concerning the sheaf of differential $1$-forms
with logarithmic singularities along $D$.
The standard references on this subject are
D\'eligne's book \cite{Del} and Saito's paper \cite{Sai}.
We define $\Omega^1_X (log D)$ as a subsheaf
of the sheaf $\Omega^1_X ( D)$ of $1$-forms with poles at most on $D$
and of order one.  This sheaf is the image of the natural map
$$ \Omega^1_X \oplus \hol_X^r
\,\,\, \longrightarrow \,\,\, \Omega^1_X (D) $$
which is given by the inclusion $\Omega^1_X \ra \Omega^1_X (D)$
and the homomorphisms sending $1 \in \hol_X \ra d log(F_i)$.
For experts we note that our definition differs from the one
in \cite{Sai} when $D$ is not normal crossing.
Saito's sheaf is the double dual of our
$\,\Omega^1_X (log D) $, which explains why
his is always locally free when $X$ is a surface
\cite[Corollary 1.7]{Sai}. Ours need not be
locally free even for surfaces.
However, our definition gives a natural exact sequence.

\begin{lemma} \label{exactlem}
If $X$ is factorial and complete then we have an exact sequence
\begin{equation}
\label{lemmasequence}
   0 \ra  \Omega^1_X \ra \Omega^1_X (log D) \ra \bigoplus_{i=1}^{r}
    \hol_{D_i} \ra 0 .
\end{equation}
\end{lemma}

\begin{proof}
The local sections of the sheaf
$\,\Omega^1_X (log D) \,$ are rational $1$-forms
which can be written as
$\,\omega = \Sigma_{i=1}^{r} \psi_i \cdot d log(F_i) +
\eta$, where $\eta$ is a regular $1$-form.

Ssince the $D_i$'s are distinct
prime divisors and $X$ is factorial, the local rings $\hol_{X, D_i}$
are discrete valuation rings with parameter $F_i$. Thus $F_j$ is invertible
in this local ring for $ j \neq i$, and $\omega$ is regular if and only
if $F_i$ divides $\psi_i$. This implies that the homomorphism
which sends $\omega$ to
the vector $\,\bigl(\psi_i \,( mod \,F_i\bigr))_{i=1,\ldots,r}\,$
is well defined, and it induces
an isomorphism from the quotient $\,\Omega^1_X (log D) / \Omega^1_X\,$
onto $\,\oplus _{i=1}^{r}   \hol_{D_i} $.
\end{proof}

Assume now that $X$ is smooth. Then both sheaves
$\Omega^1_X ( D) $ and $\Omega^1_X$ are locally free of rank $ d= dim(X)$.
Hence the intermediate sheaf $\Omega^1_X (log D) $ is  torsion free
of the same rank. Our next result shows that $\Omega^1_X (log D) $
is locally free if and only
if the divisors $D_i$ are smooth and intersect transversally.

\begin{prop} \label{Prop1}
Let $x \in X$ be a smooth point,
$x_1, \ldots, x_d$ local coordinates at $x$
and $D_1, \ldots, D_h$ the
divisors which contain $x$. Then
the sheaf $\Omega^1_X (log D) $ is
locally free at $x$ if and only if
   the $h \times d$-matrix $(\partial F_i/ \partial x_j )$
   has rank $h$ at $x$.
\end{prop}

\begin{proof}
Any local section of
$\,\Omega^1_X (log D) \,$ can be written in the form
\begin{equation}
\label{omegasum}
   \omega  \quad = \quad \sum_{i=1}^{r} \psi_i \cdot d log(F_i) +
\eta  \quad = \quad \sum_{i=1}^{h} \psi_i  \cdot d log(F_i) +
\sum_{j=1}^{d} \eta_j \cdot d x_j .
\end{equation}
This observation gives rise to a local exact sequence
\begin{equation}
\label{localexseq}
   0 \ra \hol_{X,x}^h \ra \hol_{X,x}^h \oplus \hol_{X,x}^d \ra
\Omega^1_{X,x} (log D) \ra 0.
\end{equation}
The surjective map on the right takes $((\psi_i),(\eta_j))$
to the sum on the right hand side of
(\ref{omegasum}). The injective map on the left
takes the $h$-tuple $(A_1,\ldots,A_h)$ to
$$ ((\psi_i),(\eta_j))\quad \hbox{with} \quad \psi_i = F_i A_i
\,\,\hbox{and} \,\, \eta_j = -\sum_{l=1}^h A_l \frac{\partial
F_l}{\partial x_j}. $$ The exactness of the sequence
(\ref{localexseq}) follows from the proof of Lemma \ref{exactlem}. If the
section $\omega$ in  (\ref{omegasum})
is identically zero in $\Omega^1_{X,x}(log D)$ 
then $\omega$ is in particular regular,
and so $F_i$ divides each $\psi_i$.

Now, since $X$ is reduced, a coherent sheaf $\FF$ is locally free
of rank $d$ if and only if $ dim_{\C} \FF \otimes \C_{x} = d$
for each point $x$. Since tensor product is right exact, it follows that
this condition is verified for $\Omega^1_X (log D) $  if and only if
the matrix of $\hol_X^h \ra \hol_X^h \oplus \hol_X^d$ , evaluated at
$x$, has rank precisely $h$. Since the functions $F_1, \ldots, F_h$
vanish at $x$, this is exactly the asserted condition that
the Jacobian  marix $\,(\partial F_i/ \partial x_j )_{i=1,
\dots h, j=1, \dots d}\,$  has rank $h$ at $x$.
\end{proof}

In the above situation
   where $X$ is smooth and  $\Omega^1_X (log D)$
is locally free we shall say that the divisor $D$ has
\emph{global normal crossings} (or GNC).

\begin{thm} \label{chernclassthm}
Let $X$ be smooth and assume that $D$ is a GNC divisor.
Then
\begin{enumerate}
\item
   the section $ \,d log(f)\,$ of $\Omega^1_X (log D)$
does not vanish at any point of $D$,
\item
   if the divisor $D$ intersects every curve in $X$
(in particular, if $D$ is ample) then $ d log (f)$ vanishes only
on a finite subset of $ \,V = X \backslash D$,
\item if the above conclusions hold, then
the number of critical points of $f$ on $V$,
counted with multiplicities, equals the degree of the top Chern class
$c_d (\Omega^1_X (log D)).$
\end{enumerate}
\end{thm}

\begin{proof}
We abbreviate $\,\sigma :=  dlog (f)
    = \Sigma_{i=1}^{r} u_i  \ d log (F_i)  $.
By the proof of Proposition \ref{Prop1}
   it follows that if
$(\partial F_i/ \partial x_j )_{i=1, \dots h, j=1, \dots d}\,$
has rank $h$ at $x$, then $\Omega^1_X (log D)$
is locally free of rank $d$ with generators $ d log (F_i) $ and some
choice of $d-h$ of the $d x_j$. If we write $\sigma$ in this basis,
the coefficients
of $d log (F_i)$ are the constants $u_i$ while the coefficients of
the $d x_j$ are
some regular
   functions. The first assertion follows immediately since the
exponents $u_i$ are all nonzero.

The second assertion follows from the first:
let $Z_{\sigma}$ be the zero set of the section $\sigma$.
Since $Z_{\sigma}$ does not intersect $D$, it follows that
$ dim (Z_{\sigma}) = 0$.

Thirdly, if $\FF $ is a locally free sheaf of rank $d$ on a smooth
variety $X$ of dimension $d$, and $\sigma $ is a section of $ H^0(
\FF)$ with a zero scheme $Z_{\sigma}$ of dimension $0$, then the
length of $ \,Z_{\sigma} \,$ equals the degree
of the top Chern class $c_d (\FF)$.
\end{proof}

The {\em total Chern class} of a sheaf $\FF$ is the sum
   $\, c_{tot} (\FF) =
\Sigma_{i=0}^{d}  c_i (\FF) z^i $.
This is a polynomial in $z$ whose coefficients
are elements in the Chow ring $\,A^*(X)$.
Recall that every element in $A^*(X)$ has a well-defined
{\em degree} which is the image of its
degree $d$ part under the degree homomorphism
$\,A^d(X) \rightarrow \zz$.

\begin{cor}
Suppose that $X$ is smooth and $D$ is a GNC divisor
on $X$ which intersects every curve.
Then the number of critical points of $f$,
counted with multiplicities,
is the degree of the coefficient of $z^d$
   in the following polynomial:
\begin{equation}
\label{totalChernclass}
   c_{tot} (\Omega^1_X)  \cdot \Pi_{i=1}^{r}
    ( 1 - z D_i)^{-1} \qquad \in \,\, A^*(X)[z].
\end{equation}
\end{cor}

\begin{proof}
The total Chern class
$\, c_{tot} (\FF) \,$ is multiplicative
with respect to  exact sequences, i.e.,
if $ 0  \ra A \ra B \ra C \  \ra 0$ is an exact sequence
of sheaves, then
$\,c_{tot}(B) = c_{tot}(A) \cdot c_{tot}(C)$.
Hence the sequence (\ref{lemmasequence}) implies the result.
\end{proof}

In the next section, we apply the formula
(\ref{totalChernclass}) in the case when $X$
is a smooth projective toric variety.
The Chow group
$A^d(X)$ has rank one and is generated
by the class of any point. This canonically identifies $A^d(X)$ with
$\zz$
and so any top Chern class can be considered to be a number.

\begin{cor}
\label{toricformula}
Suppose $X$ is a smooth  toric variety with boundary divisors
$\Delta_1,\ldots,\Delta_s$ and  $D$ is GNC
and meets every curve.
The number of critical points of $f$, counted with multiplicity,
equals the coefficient of $z^d$ in
\begin{equation}
\label{toricformula2}
   \frac{\Pi_{j=1}^{s}  ( 1 - z \Delta_j)}{
   \Pi_{i=1}^{r}
    ( 1 - z D_i)}
  \qquad \in \,\, A^*(X)[z].
\end{equation}
\end{cor}

\begin{proof}
By virtue of equation \eqref{totalChernclass} we need only compute
the total Chern class $\,c_{tot}(\Omega^1_X)$. For this we use the
exact sequence
in \cite[page 87]{Ful},
$$ 0 \ra  \Omega^1_X \ra \Omega^1_X (log \Delta) \ra \bigoplus _{j=1}^{s}
    \hol_{\Delta_j} \ra 0 , $$
where  $\Delta = \sum_{j=1}^s \Delta_j$, and the fact that
$\Omega^1_X (log \Delta)$ is trivial.
\end{proof}

\section{Models defined by Generic Polynomials}

We now apply the results of the previous section to  models
$\,{\bf f}  : \rr^d  \rightarrow \rr^n$. To illustrate how this works,
we first prove  Theorem \ref{thm1} for generic polynomials.
The proof of the statement that the ML degree of
generic polynomials is an upper bound on the ML degree
of special polynomials (when this number is finite)
is deferred to  Theorem \ref{toricformula3} which
is a generalization of Theorem \ref{thm1}.
See also Theorem \ref{semiconti} where this
semi-continuity principle is stated in general.

\begin{proof}[Proof of Theorem \ref{thm1} (generic case)]
The polynomials $f_1, \ldots, f_n$  are
assumed to be generic among all
(nonhomogeneous) polynomials
of degrees $b_1,\ldots,b_n$
in $\theta_1,\ldots,\theta_d$, and $u_1, \ldots, u_n$ 
are positive integers.
We take  $X$ to be projective
space $\PP^d$ with coordinates
$(\theta_0:\theta_1:\cdots :\theta_d)$.
Our object of interest is the
following rational function
on $X = \PP^d$:
$$ F \quad = \quad
   (f_1^{u_1} f_2^{u_2} \cdots f_n^{u_n})\bigl(
\frac{\theta_1}{\theta_0},
\frac{\theta_2}{\theta_0},\ldots,
\frac{\theta_d}{\theta_0}\bigr).
$$
The global factorization (\ref{globfac})
of this $F$ has $r = n+1$ prime factors,
namely,
$$ F_i \,\,\, = \,\,\,
\theta_0^{b_i} \cdot f_i(\frac{\theta_1}{\theta_0},\ldots,
\frac{\theta_d}{\theta_0}) \qquad
\hbox{for $i = 1,\ldots,n$},
$$
and $\,F_{n+1} =  \theta_0 \,$
with $\, u_{n+1} \, = \,
- b_1 u_1 - b_2 u_2 - \cdots - b_n u_n $.
The Chow ring of $X = \PP^d$ is
$\, \zz[H]/ \langle H^{d+1} \rangle$,
where $H$ represents the hyperplane class.
By our genericity hypothesis,
the $ r= n+1$ prime factors of $F$
are smooth and global normal crossing.
They correspond to the following divisor classes:
$$ D_1 = b_1 H, \,
D_2 = b_2 H,\,
\ldots,\,
D_n = b_n H \,\,\,
\hbox{and} \,\,\,
D_{n+1} = H.
$$
Projective space $\PP^d$ is a smooth
toric variety with $d+1$ torus-invariant divisors
$\Delta_j$, each having the same class $H$.
Hence the formula in (\ref{toricformula2}) specializes to
$$ \frac{(1-zH)^{d+1}}
{(1- z b_1 H) \cdots (1- z b_n H) (1- z H)}
\quad =  \quad
   \frac{(1-zH)^{d}}
{(1- z b_1 H) \cdots (1- z b_n H)}. $$
Since we work in the Chow ring of
projective space $\PP^d$, the coefficient
of $(zH)^d$ is the same as the coefficient of $z^d$
in the generating function in (\ref{thm1formula}).
\end{proof}

We now generalize our results from polynomials of fixed degrees to
Laurent polynomials with fixed Newton polytopes.  Recall that
the \emph{Newton polytope} of a Laurent
polynomial $f(\theta_1, \ldots, \theta_d)$
is the convex hull of the set of exponent vectors of the monomials appearing
in $f$ with nonzero coefficient. Given a convex polytope $P \subset \rr^d$ with
vertices in $\zz^d$, by a  \emph{generic Laurent polynomial
with Newton polytope $P$}  we will mean a sufficiently general
$\cc$-linear combination of
monomials with exponent vectors in $P \cap \zz^d$.

In the next theorem we consider
$n$ Laurent polynomials $f_1,f_2,\ldots,f_n$ having
respective Newton polytopes $P_1,P_2,\ldots,P_n$.
Because the $f_i$'s
are Laurent polynomials, i.e., their monomials may have negative
exponents, we only consider those critical points of 
 $\,f = f_1^{u_1} f_2^{u_2} \cdots f_n^{u_n}\,$
which lie in the
algebraic torus $(\cc^\ast)^d$. The number of such critical points
(counted with multiplicity) will be called the {\em toric ML degree}
of the rational function $f$.

Let $\, P \, = \, P_1 + P_2 + \cdots + P_n\,$ denote the
{\em Minkowski sum} of the given Newton polytopes,
and let $X$ be the projective toric variety defined by $P$.
Let  $\eta_1,\ldots,\eta_s \in \zz^d$ be the primitive
inner normal vectors of the facets of $P$.
They span the rays of the fan of $X$.  Let $\Delta_1,
\ldots,\Delta_s$ denote the corresponding torus-invariant divisors on
$X$. Each of the Newton polytopes $P_i$ is the
solution set of a system of linear inequalities of the specific
form
$$ P_i \quad = \quad \{\, x \in \rr^d \ \, |\ \, \langle x, \eta_j
\rangle \geq -a_{ij} \quad
\hbox{for} \quad j = 1, \ldots ,s \,\}. $$
The divisor on $X$ defined by
the Laurent polynomial $f_i$ is
linearly equivalent to $D_i = \sum_{j=1}^s a_{ij}\Delta_j$.
The $a_{ij}$ are integers which can be positive or negative.
The divisor on $X$ defined by
$f = f_1^{u_1} f_2^{u_2}\cdots f_n^{u_n}$ is linearly equivalent to
\begin{equation}
\label{divisoronX}
\sum_{i=1}^n u_i D_i \quad = \quad
\sum_{j=1}^s (\sum_{i=1}^n u_i a_{ij}) \cdot \Delta_j.
\end{equation}
We abbreviate the {\em support} of this divisor by
\begin{equation}
\label{whatisI}
I \quad = \quad
\bigl\{\, j \in \{1,\ldots,s\} \,\,|\,\,\sum_{i=1}^n u_i a_{ij} \not= 0
\,\bigr\}.
\end{equation}
A toric variety $X$ is smooth if all the cones in its normal fan are
unimodular.

\begin{thm} \label{toricformula3}
If the toric variety $X$ is smooth and
the toric ML degree of the  rational function $f$
is finite then it is bounded above by the
coefficient of $z^d$ in the following generating function with
coefficients in the Chow ring of $X$:
\begin{equation}
\label{toricchern}
   \frac{\prod_{j \notin I} (1 - z\Delta_j)}{\prod_{i=1}^n (1-zD_i)}.
   \end{equation}
Equality holds if each $f_i$
is generic with respect to its
Newton polytope~$P_i$.
\end{thm}

Note that Theorem \ref{thm1} is the special case of
Theorem \ref{toricformula3} when $P_i$
is the standard $d$-dimensional  simplex
$\,{\rm conv} \{ 0, e_1,\ldots,e_d\}\,$ scaled by
a factor of $b_i$.

\begin{proof}
Let us first assume that $f_i$ is a
generic Laurent polynomial with Newton polytope $P_i$. Let 
$\C[x_1, \dots, x_s]$ be the homogeneous coordinate ring \cite{Cox} of
$X$ with one variable for each torus-invariant divisor $\Delta_j$. Given a
Laurent polynomial $f_i(\theta)$ with Newton polytope $P_i$, the
corresponding rational function on $X$ is $F_i(x)/x^{D_i}$
where $D_i$ is as defined above and $F_i$ is homogeneous of degree
$D_i$. Therefore the rational function on $X$ we are interested in
is
$$ F \quad  = \quad  x^{-\sum u_i D_i} \prod F_i(x). $$

We next show that the divisor of $F$ is GNC. 
Note that $F_i$ is a generic 
section of a line bundle on $X$ 
that is generated by its sections. This implies,
by the Bertini-Sard theorem and by induction on $n$, that
the divisors $\{F_i=0\}$ meet transversally
in the dense torus of $X$. For points in the boundary of $X$,
we simply restrict to the torus orbit determined by the corresponding
facet where the restricted $F_i$'s remain generic sections of the
restricted bundles.

The reduced divisor of poles and zeros of $F$ is $ D = \sum D_i$ +
$\sum_{j \in I} \Delta_j$ where $I$ is defined as in (\ref{whatisI}). 
Since $\sum D_i$ is the divisor corresponding
to $P$ it is ample on $X$ by construction. So $\sum D_i$ meets
every curve on $X$ and therefore so does $D$ and we can apply
Corollary \ref{toricformula}. A variable $x_j$ appears as a factor
in $F$ if and only if $j \in I$ , in
which case $1 - z\Delta_j$ appears in both the numerator and
denominator of  (\ref{toricformula2}), and we get
the expression (\ref{toricchern}).

Consider now arbitrary Laurent polynomials $f_1, \dots, f_n$
in $\theta_1,\ldots,\theta_d$ such that
$f = \prod f_i ^{u_i}$ has only finitely many critical
points in $(\cc^*)^d$. Let $\nu$ be the coefficient
   of $z^d$ in (\ref{toricchern}).
   Let $\cc^m$ be the space of all $n$-tuples of Laurent polynomials with
the given Newton polytopes.  Consider the critical equations of $f =
\prod f_i ^{u_i}$ and clear denominators. The resulting collection of
$d$ Laurent polynomials defines an algebraic subset $\tilde{W}$ in the
product space $\cc^m \times (\cc^*)^d$. Saturate $\tilde{W}$ to remove any
components along the hypersurfaces $\{f_i = 0\}$ and
get a new algebraic subset $W$.
The map from $W$ onto $\cc^m$ is dominant and generically finite,
and the generic fiber of this map consists of $\nu$ points.

Our given Laurent polynomials $f_1, \dots, f_n$ represent a point
$\phi$ in $\cc^m$.  Let $\theta^{(1)}, \ldots, \theta^{(\kappa)}$ be
the isolated critical points of $f$. For each $i$, consider any
irreducible component $W^{(i)}$ of $W$ containing the point
$(\phi,\theta^{(i)})$ in $W \subset \cc^m \times (\cc^*)^d$.  By
Krull's Principal Ideal Theorem, the component $W^{(i)}$ of $W$ has
codimension $\leq d$ and hence it has dimension $\geq m$. As the
generic fiber is finite, the dimension of $W^{i}$ is exactly $m$ and
the projection to $\cc^m$ is dominant. Since $\theta^{(i)}$ is an
isolated solution of the critical equations, the projection map to
$\cc^m$ is open \cite[(3.10)]{Mum}, so the intersection of
$W^{(i)}$ with an open neighborhood of $(\phi,\theta^{(i)})$ maps onto
an open neighborhood of $\phi$.  Hence every generic point $\tilde
\phi$ near $\phi$ has a preimage $(\tilde \phi,\tilde{\theta^{(i)}})$
near $(\phi,\theta^{(i)})$, and these preimages are distinct for
$i=1,\ldots,\kappa$. We conclude that $\kappa \leq \nu$.
This semicontinuity argument is called  the ``specialization principle''
stated in Mumford's book  \cite[(3.26)]{Mum} and also works
when the $\theta^{(i)}$ have multiplicities, as shown in
Theorem \ref{semiconti} below. 
   \end{proof}

$\!\! $  We illustrate Theorem \ref{toricformula3} with  two
examples which we revisit in Section~5.

\begin{ex} \label{rectangles} \rm
Consider $n$ generic polynomials
$f_1(\theta_1, \theta_2), \ldots, f_n(\theta_1, \theta_2)$ where
the support of $f_i$ consists of monomials $\theta_1^p \theta_2^q$
with $0 \leq p \leq s_i$ and $0 \leq q \leq t_i$,
and suppose the $u_i$'s are generic. The Newton polytope
of $f_i$ is the rectangle
$$ P_i  \quad = \quad
\mathrm{conv}\{(0,0), (s_i, 0), (0,t_i), (s_i,t_i)\}.$$
The Minkowski sum of these rectangles is another
rectangle, and $\,X = \pp^1 \times \pp^1$.
In the numerator of (\ref{toricchern}), the contribution
of the two torus-invariant divisors $D$ and $E$ corresponding to the left
and the bottom edge of this rectangle survives.
The denominator
comes from the product of the divisors of $f_1,\ldots,f_n$:
$$ \frac{(1-zD)(1-zE)}{(1-(s_1D+t_1E)z)(1-(s_2D+t_2E)z) \cdots
(1-(s_nD+t_nE)z)}.$$
Now, the coefficient of the term $z^2$ modulo the Chow ring relations
$$ D^2 = 0, \quad E^2 = 0, \quad D \cdot E = 1$$
gives the toric ML degree
\begin{equation}
\label{rectangleformula}
  (\sum_{i=1}^n s_i)(\sum_{j=1}^n t_j) \, + \, \sum_{k=1}^n s_k t_k
\,\, - \,\,
   \sum_{i=1}^n (s_i+t_i) \,+\, 1.
\end{equation}
\end{ex}

\begin{ex}
\label{haseightrays}
  Let $f_1, f_2, f_3$  be generic polynomials
in $\theta_1$ and $\theta_2$ with supports
\begin{align*}
A_1 & \quad = \quad \{1, \theta_1, \theta_1\theta_2, \theta_1^2 \} , \\
A_2 & \quad = \quad \{1, \theta_1, \theta_2, \theta_1\theta_2,
\theta_1^2\} , \\A_3 & \quad = \quad \{1, \theta_1\theta_2,
\theta_1\theta_2^2\} .
\end{align*}
The corresponding Newton polytopes $P_1,P_2,P_3$
are shown in Figure \ref{fig2}. \begin{figure}[h]
\centerline{
\epsfig{file=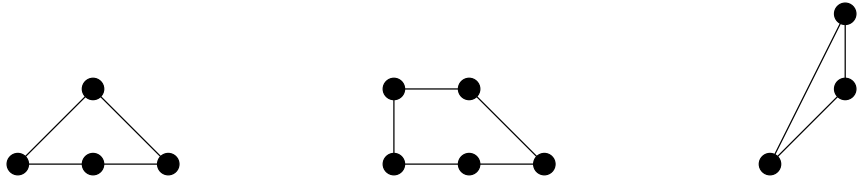}}
\caption{Three Newton polygons}
\label{fig2}
\end{figure}

\noindent
The  normal fan of the Minkowski sum has eight rays and is shown
in Figure~\ref{fig3}.
\begin{figure}[h]
\centerline{
\epsfig{file=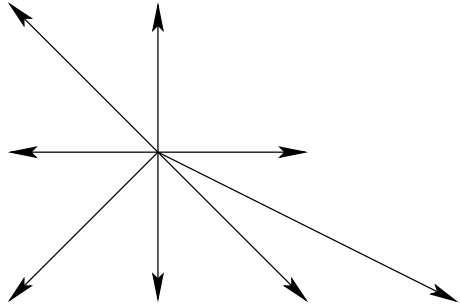}}
\caption{The fan of a smooth projective toric surface}
\label{fig3}
\end{figure}

Theorem \ref{toricformula3} applies because the toric surface $X$ is smooth.
We label the eight rays by $x_1,\ldots,x_8$ in counterclockwise order,
starting with $(1,0)$.
The Chow ring $\,A^*(X) \,$ is the polynomial ring
$\zz[x_1,\ldots,x_8]$ modulo the ideal
\begin{align*}
\langle &\, x_1x_3, \,x_1x_4, \,x_1x_5, \, x_1x_6,\, x_1x_7,
\,x_2x_4, x_2x_5, \,x_2x_6,
\, x_2x_7, \,x_2x_8,
\\ & x_3x_5,\, x_3x_6, \,x_3x_7, \, x_3x_8,\, x_4x_6, \,x_4x_7,\, x_4x_8,
\,x_5x_7, \,x_5x_8, \,x_6x_8,
   \\ & x_1 - x_3 - x_4 - x_5 + x_7 + 2x_8 \,,\, \, x_2 + x_3 - x_5 -
x_6 - x_7 -x_8\, \rangle.
\end{align*}
The three divisors corresponding to the polygons
$P_1, P_2, P_3$ in Figure \ref{fig2} are
\begin{align*}
D_1 &\quad = \quad 2x_3 + 2x_4 + 2x_5 + x_6 \\
D_2 &\quad = \quad 2x_3 + 2x_4 + 2x_5 + x_6 + x_7 + x_8 \\
D_3 &\quad = \quad x_4 + 3x_5 + 2x_6 + x_7
\end{align*}
If all $u_i$ are positive, then
   the support of the divisor
$u_1 D_1 + u_2 D_2 + u_3 D_3$ is
   $I = \{3, \ldots, 8\}$. It follows that the
toric ML degree is the coefficient of
$z^2$ in
$$ (1-zx_1)(1-zx_2)(1-zD_1)^{-1}(1-zD_2)^{-1}(1-zD_3)^{-1}. $$ This
coefficient is $\, 14 x_1 x_2 $, which means that the
toric ML degree is $14$.  \qed
\end{ex}

The {\em toric ML degree of the model} $\,{\bf f}\,$
is the toric ML degree defined above for generic $u$.
In this case, there is no cancellation among the
coefficients in (\ref{whatisI}), and $I$ is the set of
all indices $j$ such that for some $P_i$ the supporting hyperplane
normal to $\eta_j$ does not pass through the origin.
The toric ML degree of ${\bf f}$ is a
numerical invariant of the polytopes $P_1, \ldots, P_n$.
A combinatorial formula for this invariant
will be presented in Theorem \ref{toricformula4}
of Section 5.

\section{Bounded Regions in Arrangements}

As in the Introduction, we consider $n$ polynomials $f_1,\ldots, f_n$
in $d$ unknowns $\theta_1, \ldots, \theta_d$. We now assume that all
coefficients of the $f_i$'s are real numbers, and we also assume that
$u_1, \ldots, u_n$ are positive integers. However, we do not assume
that the union of the divisors of the $f_i$'s has global normal crossings. 
This is the case of interest in statistics. Consider the arrangement of
hypersurfaces defined by the $f_i$'s and let $V_{\rr} = \rr^d \setminus
\bigcup_{i=1}^n \{f_i = 0\}$ be the complement of this arrangement. A
connected component of $V_{\rr}$ is a {\em bounded region} if it is
bounded as a subset of $\rr^d$. Then the following observation
holds.

\begin{prop} \label{bounded}
For any polynomial map ${\bf f} : \rr^d \rightarrow \rr^n$ and any $u
\in \nn_{>0}^n$,
$$
\begin{array}{c} \#\{\mbox{bounded regions of } V_{\rr} \} \\
\leq \quad \#\{\mbox{critical points of }  f_1^{u_1} \cdots f_n^{u_n}
\mbox{ in }
\rr^d \} \\
\leq \quad \mbox{ ML degree of  } {\bf f}.
\end{array}$$
\end{prop}
\begin{proof} The function $f = f_1^{u_1} \cdots f_n^{u_n}$ is continuous,
and on the boundary
of the closure of each bounded region its value is zero. Hence
it has to have at least one (real) critical point in the
interior of each region. The second inequality holds trivially,
since the ML degree was defined as the number of critical 
points of $\, f_1^{u_1} \cdots f_n^{u_n}\,$ in
$\cc^d$, counted with multiplicities.
\end{proof}

This observation raises the question whether the inequalities above could be
realized as equalities. We next show that this is the case when
   $f_1, \ldots, f_n$ are quadrics in the plane. Here the
   ML degree is $2n^2-2n+1$ by Theorem \ref{thm1}.

\begin{prop}For each $n$, there are $n$ quadrics $f_1, \ldots, f_n$
in $\rr^2$
such that
$$ \#\{\mbox{bounded regions of } V_{\rr} \} \quad
= \quad \mbox{ ML degree of  } {\bf f}
\quad = \quad 2n^2 - 2n + 1. $$
Hence all critical points are real.
\end{prop}

\begin{figure} 
\centerline{
\epsfig{file=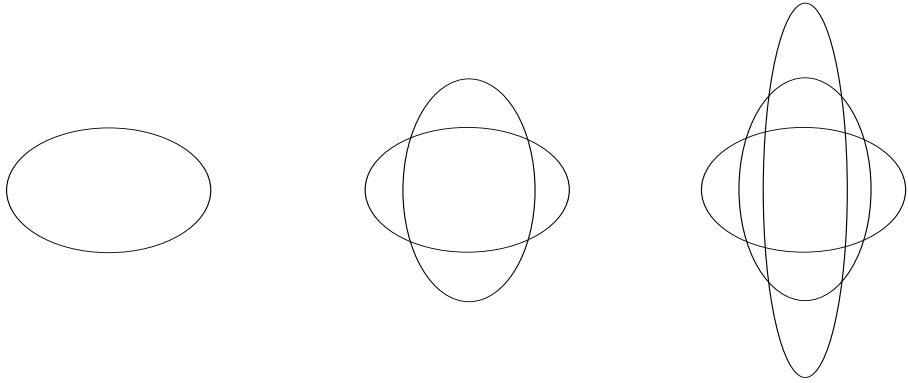, height={5cm}, width={12cm}}}
\caption{The ``nested'' ellipse construction}
\label{fig1}
\end{figure}

\begin{proof} We will take $n$ quadrics that define ``nested''
ellipses with center at the origin, as suggested by Figure \ref{fig1}.
The proof follows by induction: assume we have $2(n-1)^2 - 2(n-1) + 1$ bounded
regions with $n-1$ ellipses. Observe that the $(n-1)$st ellipse contains
$2n-3$ bounded regions. Then we add a new long and skinny ellipse
which replaces the $2n-3$ regions with $3(2n-3)+2$ regions. The
total count comes out to be $2n^2 - 2n + 1$.
\end{proof}

We will see such an equality holding for $n$ linear hyperplanes in
$\rr^d$ below. However, even in the plane  $\rr^2$, the number of critical
points and the number of bounded regions
of $V_{\rr}$ diverge for curves of degree $\geq 3$.
Theorem \ref{thm1} implies that for $n$ generic plane curves of
degrees $b_1, \ldots, b_n$ the ML degree is
$$\sum_{i=1}^n b_i(b_i-2) \, + \, \sum_{i<j}b_ib_j  \, + \, 1.$$
The optimal upper bound for the number of bounded regions of $V_{\rr}$
is smaller than the ML degree, by the following unpublished
result due to Oleg Viro.

\begin{thm} \label{thm:viro}
{\rm (Viro) }
Let $f_1, \ldots, f_n$ be real plane curves
of degrees $b_1, \ldots, b_n$, and let $K$ be the number of
odd degree curves among them. Then
$$\#\{\mbox{bounded regions of } V_\rr\} \quad \leq  \quad
\sum_{i=1}^n \frac{(b_i-1)(b_i-2)}{2} \, + \, \sum_{i<j}b_ib_j
\, + \, 1 - K,$$
and this bound is optimal.
\end{thm}

\begin{proof} The proof is by induction. For $n=1$ the
bound above is Harnack's inequality \cite{Har}, and it is optimal.
Suppose the statement is proved for $n-1$ curves, and $f_n$ defines
a curve of degree $b_n$. We will take this new curve so that it has
the maximum number of bounded regions allowed by
Harnack's inequality,
i.e. $B_n := (b_n - 1)(b_n-2)/2  + 1$ if $b_n$ is even, and one
less than that if $b_n$ is odd. One can achieve $B_n$ by taking
a curve with $B_n-1$ unnested ovals and one more distinguished piece
(that gives an extra bounded region when $b_n$ is even) such that
some line intersects this distinguished piece in exactly $b_n$ points.
We can arrange this last curve in such a way that when
we superimpose the distinguished piece on the arrangement
given by $n-1$ curves, the last curve will intersect the curve
given by $f_i$ in $b_ib_n$ points. Now if we trace this
last distinguished piece, every time we encounter an intersection
point an extra bounded region is created, except the last point in case
$b_n$ is odd. Together with the remaining $B_n-1$ ovals we get
the bound.
\end{proof}

In order to get any meaningful lower bound on the number of bounded
regions of $V_\rr$ one needs to make some assumptions. Without
any assumptions the lower bound is zero: for $f_i$ of even degree
we take an empty (real) curve, and for $f_i$ of odd degree we take the
union of an empty curve with a line. If we let all the lines intersect
in a single point there will not be any bounded region. If we insist
on at least having a GNC configuration, then by the same construction
the lower bound we get is the number of bounded regions
in a generic arrangement of $K$ lines where $K$ is the
number of odd degrees $b_i$.
This idea leads us to studying
the ML degree of a hyperplane arrangement.

\begin{thm} \label{thm-linear} Let ${\bf f}$ be given by
$n$ linear polynomials $f_1,\ldots,f_n$ with real coefficients.
Then the ML degree of ${\bf f}$  is equal to the number of
the bounded regions of $V_\rr$, and all critical points of
the optimization problem (\ref{mle}) are real.
\end{thm}

This theorem does not assume any hypothesis such as
global normal crossing. Under the GNC hypothesis,
the hyperplanes would be in general position and
the number of bounded regions equals $\binom{ n-1}{  d}$,
as predicted by Theorem \ref{thm1}.
Theorem \ref{thm-linear} is essentially due to
   Varchenko \cite{V}. We shall give a new proof.

\begin{proof}
In light of Proposition \ref{bounded},
we need to show that the number of bounded regions
of $V_\rr$ equals the number of complex solutions of the
ML equations.  Let $\,f_i = \sum_{j=1}^d a_{ij}\theta_j + c_i\,$
for $i=1,\ldots,n$.
The ML equations are
\begin{equation} \label{new-eq}
\sum_{i=1}^n \frac{u_i a_{i1}}{f_i} \, = \, \cdots \,
= \, \sum_{i=1}^n \frac{u_i a_{id}}{f_i} \,\, =\,\, 0.
\end{equation}
Consider the map $\psi : \cc^{d+1}  \rightarrow \cc^n$
given by $\psi(\theta_0, \ldots, \theta_d) =
(1/F_1, \ldots, 1/F_n)$.
Here $F_i = c_i\theta_0 + \sum_{j=1}^d a_{ij}\theta_j\,$ is
the homogenization of $f_i$. We let $\bar{\mathcal H}$
be the central hyperplane arrangement in $\rr^{d+1}$ given
by the $F_i$. We assume that the intersection of all
the hyperplanes in $\bar{\mathcal H}$ is just the origin;
otherwise, the linear forms $F_i$ depend on fewer than $d$ 
coordinates, and then we get infinitely many critical points.
The Zariski closure of $\mathrm{im}(\psi)$ in $\pp^{n-1}$
is a $d$-dimensional complex variety $\mathcal{V}$. The solution set
on $\mathcal{V}$ of the $d$ linear equations
$$\sum_{i=1}^n (u_i a_{i1}) y_i \, = \, \cdots \, = \,
\sum_{i=1}^n (u_i a_{id}) y_i \,\,=\,\, 0$$
consists of finitely many  points
provided $u_1, \ldots, u_n$ are generic.
   Obviously, the solutions to (\ref{new-eq}) lift
to such complex solutions. In other words, the degree of
the projective variety $\mathcal{V}$ is
an upper bound on the ML degree of ${\bf f}$.

   Now we will compute  the degree of $\mathcal{V}$.
This variety is the projective spectrum of the
$\nn$-graded algebra $\,R =  \cc[1/F_i:
\, i=1, \ldots, n]\,$ where $\mathrm{deg}(1/F_i) = 1$.
Terao  \cite[Theorem 1.4]{T} showed that
the Hilbert series of $R$ is equal to
\begin{equation}
\label{Terao}
\sum_{X \in \mathcal{L}} (-1)^{\mathrm{codim}(X)} \mu(X)
\left( \frac{t}{1-t} \right)^{\mathrm{codim}(X)}
\end{equation}
Here $\mu$ is the M\"obius function of the intersection 
lattice $\mathcal{L}$ of the  arrangement $\bar{\mathcal{H}}$. From 
this series we shall determine the  leading coefficient of its
Hilbert polynomial. This coefficient has
the form $e/d!$, where $e$ is the degree of $\mathcal{V}$.

For large enough $r$, the coefficient of $t^r$ in the
Hilbert series  (\ref{Terao}) equals
\begin{equation}
\sum_{i=0}^{d+1} (-1)^i \binom{r-1}{i-1}
\sum_{\mathrm{codim}(X) = i} \mu(X).
\end{equation}
This is the Hilbert polynomial
of the graded algebra $R$. Its leading term is
$$(-1)^{d+1} \mu(0) \frac{r^d}{d!}.$$
We conclude that the degree of
the projective variety $\mathcal{V}$
is $(-1)^{d+1} \mu(0)$.
By Zaslavsky  \cite{Zas},
this number equals the number of
bounded regions of $V_\rr$.
\end{proof}

\begin{ex}
A family of important statistical models where
Theorem \ref{thm-linear} applies
is the {\em linear polynomial model} of \cite{PRW}.
Such a model is given by a polynomial
in $r$ unknowns $x = (x_1,\ldots,x_r)$ with indeterminate
coefficients,
$$ p(x) \quad = \quad \sum_{j=1}^d \theta_j x^{a_{j}}
\qquad \hbox{(with $a_j \in \nn^r$)}, $$
together with $n$ data points
$v_1, \ldots, v_n \in \rr^r$.
The model is parametrized by
$$ f_1 (\theta) \, = \,
\sum_{j=1}^d  \theta_j v_1^{a_j}, \,\,\ldots, \,\,
   f_n (\theta) \, = \,
\sum_{j=1}^d  \theta_j v_n^{a_j}. $$
The ML degree is the number of bounded regions
of this arrangement.  \qed
\end{ex}

\section{Polytopes and Resolution of Singularities}

We now return to the setting of Section 3, with the aim of
relaxing the restrictive smoothness hypothesis in Theorem
\ref{toricformula3}. Our aim is to derive a combinatorial formula
for the toric ML degree of any model $\, {\bf f} \,$ defined by
generic Laurent polynomials satisfying a mild hypothesis. The
derivation of Theorem \ref{toricformula4} involves resolution of
singularities in the toric category. In the end of the section we
shall comment on using resolution of singularities for bounding
the ML degree in general.

Given a polytope $P$ in $\rr^d$ and a linear functional $v$ on
$\rr^d$, we write
$$ P^v \quad = \quad \bigl\{ \,p
\in P \,\,|\, \
\forall p' \in P \,:\,
\langle v, p \rangle \leq \langle v, p'\rangle \  \bigr\} $$
for the face of $P$ at which $v$ attains its minimum.
Two linear functionals $v$ and $v'$ are
{\em equivalent } if $\,P^v = P^{v'}$. The equivalence
classes are the relative interiors of cones of the
{\em inner normal fan} $\Sigma_P$.
If $\sigma$ is a cone in  $\Sigma_P$,
or $\sigma$ is a cone in any fan which refines $\Sigma_P$,
then we write $P^\sigma = P^v$ for $v $ in the
relative interior of $\sigma$.
If $f$ is a polynomial
with Newton polytope $P$
then $f^\sigma$ denotes the {\it leading form}
consisting of all terms of $f$ which are supported
on $P^\sigma$.

As in Section 3, let $f_1,\ldots,f_n$ be Laurent polynomials
with Newton polytopes $P_1,\ldots,P_n \subset \rr^d$. Consider
any fan $\Sigma$ which is a common refinement of the inner normal fans
$\Sigma_{P_1},\ldots,\Sigma_{P_n}$. Suppose
$\tau$ is a cone in $\Sigma$ and let
$k$ be the dimension of $(P_1+\cdots+P_n)^{\tau} $.  There exists a
$k$-dimensional linear subspace $L$ of $\rr^d$
and vectors $q_1,\ldots,q_n \in \rr^d$ such that
$\,q_i + P_i^\tau \,$ lies in $L$ for all $i = 1,\ldots,n$.
The subspace $L$ is unique and satisfies
$L \cap \zz^d \simeq \zz^k$. Let
$\,V(\,\cdot \, ,\ldots,\,\cdot \,)\,$ denote
the {\em normalized  mixed volume} on the subspace $L$.
Here ``normalized'' refers to the lattice
$L \cap \zz^d$, as
is customary in toric geometry \cite{Ful}.
For any $k$-element subset $\{i_1,\ldots,i_k\}$ of
$\{1,2,\ldots,n\}$ we abbreviate
\begin{equation}
\label{mixedvolume}
  V(P_{i_1},\dots, P_{i_k}; \tau)
\,\,=\,\, V(q_{i_1} + P_{i_1}^{\tau}, \dots,
q_{i_k} +  P_{i_k}^{\tau}) \,
\quad \hbox{if $\,\rm{codim}( \tau) = k$}, 
\end{equation}
and $\, V(P_{i_1},\dots, P_{i_k}; \tau)
= 0 \,$ if $\, \rm{codim}(\tau) > k$.
If $k = d$ and $\tau = \{0\}$ then
we simply write $\, V(P_{i_1},\dots, P_{i_d})\,$
for the mixed volume in  (\ref{mixedvolume}).
If $k = 0$ and $\tau $ is full-dimensional
then (\ref{mixedvolume}) equals $1$; this happens in
the last sum of (\ref{longmultline}).

We are now ready to state our more general toric ML degree
formula. As in Section 3, let $X$ be the toric variety
corresponding to the Minkowski sum $P = P_1 + \dots + P_n$ and
$\Sigma_X$ the normal fan with rays $\eta_1, \dots, \eta_s$.
We consider the function $f = f_1^{u_1}
\cdots f_n^{u_n}$. Each polytope $P_i$ corresponds to a divisor
$D_i$ so the divisor of $f$ is $D = \sum u_i D_i$. Let $I$ be the
support of $D$ as in \eqref{whatisI}. Label the rays of $\Sigma_X$
so that $\{1, \dots, r\}$ are the indices not in $I$.

For each subset $J$ of $\{1,\ldots,r\}$ let $\tau_J$ denote the
smallest cone of $\Sigma$ which contains the vectors $\,\eta_j \,$
for $j \in J$. If no such cone exists then $\tau_J$ is just a
formal symbol and  the expression (\ref{mixedvolume})  is declared
to be zero for $\tau = \tau_J$. The mild smoothness hypothesis we
need is that every singular cone of $\Sigma$ contains at least one
ray from $I$. Equivalently all cones $\tau_J$ are smooth.

  \begin{thm} \label{toricformula4}
Suppose every singular cone of $\Sigma_X$ contains some ray in the
support of the divisor $D$. Then, the toric ML degree of the rational 
function $f$ is bounded above
by the  following  alternating sum of mixed volumes:
\begin{multline}
\label{longmultline}
  \sum_{1 \leq i_1 \leq \cdots \leq i_d \leq n} V(P_{i_1}, \dots,
P_{i_d}) \,\,\, - \!
\mathop{\sum_{j \in \{1, \dots, r\}}}_{1 \leq i_1 \leq \cdots
\leq i_{d-1} \leq n} V(P_{i_1}, \dots, P_{i_{d-1}}; \tau_{\{j\}})  \\
+ \,\,\mathop{\sum_{\{j_1,j_2\} \subset \{1, \dots, r\}}}_{1 \leq i_1 \leq
\cdots \leq i_{d-2} \leq n} V(P_{i_1}, \dots, P_{i_{d-2}};
\tau_{\{j_1,j_2\}}) \quad + \quad \cdots \quad + \\
(-1)^d \sum_{\{j_1, \dots, j_d\} \subset \{1, \dots, r\}} V(\emptyset;
\tau_{\{j_1,\ldots,j_d\}}).
\end{multline}
Equality holds if each $f_i$ is generic relative
to its Newton polytope~$P_i$.
\end{thm}

\begin{proof}
In order to apply Corollary \ref{toricformula} we must resolve the
singularities of $X$. For toric varieties this is done in two
steps. First we get a simplicial toric variety without adding any new
rays to the fan.  Second we resolve the remaining singular (but
simplicial) cones by adding new rays. This procedure is described
in detail in \cite{Ful}. Typically the first step involves
taking the pulling subdivision at each ray in the fan. However,
under the given hypothesis it is enough to perform pulling
subdivisions only at the rays in the support of $D$ to obtain
a simplicial fan $\Sigma_{\tilde{X}}$. This fine detail
will be important below. Our hypothesis
holds for this intermediate fan as well, and subsequently
we take a smooth refinement $\Sigma_{X'}$ of
$\Sigma_{\tilde{X}}$ by adding new rays in the relative interiors
of each of the singular cones. Let $\pi \,: \, X' \to X$ be the
induced map.

We will show that we get no new critical
points under the resolution. Hence the number of critical points
can be computed on $X'$. We finally claim that the Chern
class formula expands into the given combinatorial formula.

We investigate critical points of the pullback of our rational
function:
$$F' \,\,= \,\,\pi^{\ast}(F) \quad = \quad (x^{-\sum u_i 
\pi^{\ast}(D_i)}) \prod
\pi^\ast(F_i(x)). $$
For generic $f_i$, the same argument as in the proof of
Theorem \ref{toricformula3} shows that the reduced
divisor of poles and zeros $D'$ of $F'$ is GNC. What we
must show is that
all critical points of $F'$ on $X'$ are off the exceptional locus hence
actually critical points of $F$ on $X$.

There are two types of new cones in $\Sigma_{X'}$.
The first come from the triangulation step. By our
construction, any such cone must contain a ray $\eta_j$ in the
support of $D$. This ray corresponding to
the strict transform under $\pi$ is in the support of
$\sum u_i \pi^{\ast}(D_i)$, and its variable appears as a factor in
$F'$. By part 1 of Theorem \ref{chernclassthm}, $F'$ has no
critical points along the torus-invariant divisor
$\Delta_j'$, hence no critical points on any torus orbit contained
in $\Delta_j'$.

The second type of new cone comes from the desingularization
step. These cones all contain at least one new ray $\eta_E$
corresponding to an exceptional divisor $\Delta_E$ in $X'$. We will
show there are no critical points on $\Delta_E$. Equivalently we show
that there are no critical points on each torus orbit contained in
$\Delta_E$.

Given a torus orbit let $\tau_E$ be the corresponding cone of
$\Sigma'$ containing $\eta_E$. There is some minimal
cone $\tau$ of $\Sigma$ containing $\tau_E$. Let $\tau'$ be any cone
of $\Sigma'$ containing $\tau_E$ that is maximal with respect to being
contained in $\tau$. Since $\tau$ is refined in $\Sigma'$ it must be a
singular cone, and so by the hypothesis it has some generating ray
in the support of $\sum u_i D_i$, or equivalently the linear function
of this Cartier divisor is not identically zero on $\tau$. The
pullback keeps the same linear functional which cannot be zero on the
subset $\tau'$ of $\tau$.  As a consequence, $\tau'$ contains a
ray $\eta_j$ in the support of $\sum u_i \pi^{\ast}(D_i)$. 
This means $x_j^c$ appears as a factor in $F'$ for some nonzero 
integer $c$. 

If $\eta_j$
is a generator of $\tau_E$, then as above there are no critical
points on $\Delta_j'$ and thus no critical points on the orbit
corresponding to $\tau_E$.
Suppose on the other hand $\eta_j$ is not a generator of $\tau_E$.
Let $x_{E_1}, \dots, x_{E_k}$ be the variables corresponding to the
generators of $\tau_E$ in an affine chart of $X'$ that contains
$\tau'$. Note that the variable $x_j$ corresponding to
$\eta_j$ is not among these variables.  
Because $\tau$ is the minimal cone of $\Sigma$ containing
$\tau_E$, the face $P_i^{\tau_E}$
is contained in the face $P_i^{\tau}$, and hence it is contained in
the face $P_i^{\tau'}$. So $(F'_i)^{\tau_E}$, obtained by setting 
$x_{E_1}, \dots, x_{E_k}$ to
zero, does not contain any of the variables corresponding to $\tau'$;
in particular it doesn't contain $x_j$. On the other hand,
$F'_i = (F'_i)^{\tau_E}+ G'_i $ where $G'_i$ is in the ideal generated by
$x_{E_1}, \dots, x_{E_k}$. Since $x_j$ is not among the $x_{E_i}$ we have
$${\bigl(\frac{\partial F'_i}{\partial x_j}\bigr)}|_{x_E = 0} = \frac{\partial
((F_i')^{\tau_E})}{\partial x_j} = 0,$$
where $x_E = 0$ means $x_{E_1} = \cdots = x_{E_k} = 0$. 
Hence we have
$$(\frac{\partial(log F')}{\partial x_j})|_{x_E=0} \quad = \quad
\frac{c}{x_j} + \sum_i {\bigl(\frac{\partial F'_i}{\partial 
x_j}\bigr)}|_{x_E = 0} \quad \neq \quad 0. $$
We conclude that $F'$ has no critical points on the torus orbit
corresponding to $\tau_E$ as desired.

Thus the toric ML degree of $f$ on $X'$ is the same as that on $f$.
Since $F'$ has no critical points on the exceptional locus and $D$
is ample on $X$, $\pi^{\ast}(D)$ meets any curve off the exceptional
locus and therefore the ML degree must be finite. It is computed
in the Chow ring of $X'$ as the coefficient of $z^d$ in
$$ \frac{\prod_{j \not \in I} (1 - z\Delta'_j)\prod_k (1-
z\Delta_{E_k})}{\prod_{i=1}^n (1 - zD'_i)}.$$
  Here $\Delta'_j$ are the strict transforms of the $\Delta_j$
not in the support of $\sum u_i D_i$. The $\Delta_{E_k}$ are
those exceptional divisors which are not in the support of $\sum u_i
\pi_i^{\ast}(D_i)$, and $D'_i$ are the proper transforms of the
divisor classes of the~$F_i$'s.

We can now expand our Chern class product replacing
$(1-zD'_i)$ in the denominator by $1 + zD'_i + z^2(D'_i)^2 + \cdots +
z^d(D'_i)^d$ in the numerator. The intersection product of any collection
of prime torus-invariant divisors is the cycle of the cone they span or
$0$ if there is no such cone. Hence the coefficient of $z^d$ is the
sum of all terms of the form
$$ (-1)^k D'_{i_1}D'_{i_2} \cdots D'_{i_{d-c}} \cdot \tau^c. $$
Here $1 \leq i_1 \leq \cdots \leq i_{d-c} \leq n$ and $\tau^c$
ranges over all dimension $c$ cones of $\Sigma'$ spanned
by rays not in the support of $\sum u_i D'_i$. This product
is exactly the mixed volume $V(P_{i_1}, \dots, P_{i_{d-c}}; \tau^c)$.

To finish we note that if $\tau^c$ contains an exceptional divisor
$\Delta_{E_k}$, the minimal cone $\tau$ of $\Sigma$ containing $\tau^c$ must
have dimension strictly larger than $c$. This is because $\tau^c$ does
not have any rays in the support of $\sum u_i D'_i$ hence is not
maximal in $\tau$. As a consequence all of the faces
$P_{i_j}^{\tau^c}$  have a translate that lies
in a subspace of dimension strictly
less than $d-c$ and the corresponding mixed volume is $0$. In conclusion,
the exceptional divisors do not contribute to the top Chern class product
and the formula reduces to the stated one.
\end{proof}

In two variables we recover a particularly simple formula:
\begin{cor} \label{areaarea}
Let $f_1 , \dots, f_n$ be generic Laurent polynomials in two
variables $(\theta_1,\theta_2)$ with Newton polygons $P_1, \dots,
P_n$. If the origin lies on none of the lines spanned by edges of 
their Minkowski sum $P = P_1 + \dots + P_n$ then
the toric ML degree equals the area of $P$ plus the areas
of each of the $P_i$.
\end{cor}

\begin{proof}
This is a special case of Theorem \ref{toricformula4} when no facets pass
through the origin. Therefore the only term is
$$\sum_{i=1}^n \sum_{j=i}^n V(P_i, P_j).$$
The Euclidean area of each $P_i$ is $V(P_i, P_i)/2$. The Euclidean area of
the Minkowski sum $P = \sum P_i $ equals
$\,\frac{1}{2}\sum_i V(P_i, P_i) + \sum_{i<j} V(P_i, P_j)$.
  The stated formula is the sum of these expressions.
\end{proof}

Now that we are equipped with  the volume
formulas in Theorem \ref{toricformula4}
and Corollary \ref{areaarea}, let us
revisit the two-dimensional examples from Section 3.

\begin{ex}
The Newton polygons $P_1,\ldots \!, \! P_n$ in Example~\ref{rectangles}
are axis-parallel parallelograms.
The first term in the formula (\ref{rectangleformula})
is the area of their Minkowski sum $\,P_1 + \cdots + P_n$,
and the second term is the sum of the areas of the $P_i$, as in
Corollary \ref{areaarea}. The third and fourth term
are the two correction terms stemming from the
fact that the origin is a vertex of each Newton polytope.
These terms disappear if we replace
one $f_i$ by $\theta_1 \theta_2 f_i$.

The number $14$ in Example \ref{haseightrays}
can also be derived  using Theorem
\ref{toricformula4}. The three polygons in
Figure \ref{fig2} have
areas $1, \frac{3}{2},$ and $\frac{1}{2}$ respectively. The area of
their sum is $15$. The two divisors corresponding to $x_1$ and $x_2$
pass through the origin and yield correction terms of $1$ and $4$
respectively. Finally we add back $1$ for the vertex at the
origin. Altogether $1 + \frac{3}{2} + \frac{1}{2} + 15 - 1 - 4 + 1 =
14$. \qed
\end{ex}

\vskip .1cm

Our discussion so far indicates that we get the sharpest results
when $X$ is smooth and $D$ is GNC. In the toric
case the smoothness hypothesis could  be largely removed  
as in Theorem \ref{toricformula4}. The GNC condition can also be relaxed
for certain other cases as we saw in the previous section. In general,
if the pair $(X,D)$ does not satisfy the smoothness and
GNC hypotheses then we must appeal to
Hironaka's theorem on resolution of singularities
(see e.g.~\cite{Hau}).
This furnishes a proper
projective morphism $ \pi : X' \ra X$ such that
$X'$ is smooth and $ \pi^{-1} (D)$  has  GNC.
We need to compute   the divisor
$\pi^* (div(f)) $ of the pullback of the function $f$.
If $D'_i$ is the proper transform of the divisor $D_i$
and $E_1,\ldots,E_k$ are the exceptional divisors of
$\pi$ then
$$ \pi^* (div(f)) \quad = \quad \sum_{i=1}^{r} u_i D'_i \,+\,
\sum_{j=1}^{k} \mu_j E_j, $$
where $\mu_j$ are certain (possibly negative)  integers.
These integers are $\mu_j = \Sigma_{i=1}^r u_i m_{ij}$
where $ m_{ij}$ is the 
multiplicity of the full transform of $D_i$ along $E_j$.
The underlying reduced divisor is
$\, D' \, :=\,  \sum_{i=1}^{r}  D'_i \,+\,
   \sum_{j: \mu_j \neq 0} E_j$.
The number of critical points is now
gotten by applying Theorem \ref{chernclassthm} to $(X',D')$
instead of $(X,D)$.  This procedure can be very complicated
in practice. We illustrate it with a simple example.

\begin{ex} \label{4 curves} Let $d=2$, $n=4$,
$\, f_1 = x$, $f_2 = y$, $f_3 = (x-1)^2 + (y-1)^2 -2$, and $f_4 =
(x+1)^2 + 2(y-2)^2 - 9$. The divisor $D$ is not a GNC divisor since at
the origin all the four curves defined by $f_1, \ldots, f_4$ meet. In
order to resolve this singularity we blow up $X = \PP^2$ at $(0:0:1)$
to obtain $X'$ which is smooth. We note that $c_{tot}(\Omega_{X'}^1)$
is $(1-zE)(1-zH)(1-zH')^2$
where $H$ is the proper transform of the generic
hyperplane section, $E$ is the exceptional divisor, and $H'$ is the
proper transform of a line through the origin.

We have four cases: we consider first the general case  where
   $u_1+u_2+u_3+u_4 \neq 0$ and  $u_1+u_2+2 u_3+ 2u_4 \neq  0$.
In this case  $D'$ consists of the proper transforms of
the four original curves, the exceptional divisor,
and the pullback of the line at infinity.
After cancellations, we
just need to compute the coefficient of  $z^2$ in
$\frac{1}{(1-C_1z)(1-C_2z)}$
   where $C_1$ and $C_2$ are the irreducible divisors
corresponding respectively to the circle and the ellipse.  This coefficient is
$C_1^2 + C_2^2 + C_1 \cdot C_2 $. In $X'$ the two
curves intersect in three points, and their self-intersection also yields three
points. Hence the ML degree is nine.

In the special case where $u_1+u_2+u_3+u_4 = 0$,
we need  to compute the coefficient of  $z^2$ in
$\frac{1-zE}{(1-C_1z)(1-C_2z)}$, which is
$C_1^2 + C_2^2 + C_1 \cdot C_2 - E \cdot C_1 - E \cdot C_2$.
Since $E \cdot C_i = 1 $  the ML degree drops down to seven.
If $u_1+u_2+2u_3+2u_4 = 0$ then
the coefficient of  $z^2$ in
$\frac{1-zH}{(1-C_1z)(1-C_2z)}$ is
$C_1^2 + C_2^2 + C_1 \cdot C_2 - H \cdot C_1 - H \cdot C_2$.
Since $H \cdot C_i = 2 $, the ML degree is five.
Finally, if both $u_1 + u_2 + u_3 + u_4 $ and $u_1 + u_2 + 2u_3 + 2u_4 $
are zero, then the coefficient of $z^2$ in 
$\frac{(1-zH)(1-zE)}{(1-C_1z)(1-C_2z)}$ is
$C_1^2 + C_2^2 + C_1\cdot C_2 - H\cdot C_1 - H\cdot C_2 - E \cdot C_1 
- E \cdot C_2 + H \cdot E$: since $H \cdot E = 0$, the ML degree 
further drops
down to three.

The number of bounded
regions of the complement of the four curves in $\rr^2$ is seven.
By Proposition \ref{bounded}, this is a lower bound on
the ML degree when all $u_i$ are positive. This example shows that,
for specific {\em negative} values of the $u_i$'s, the number of critical
points may be smaller than this lower bound.
\qed
\end{ex}

\section{ML degree and Euler characteristic}

A well-known result in the theory of hyperplane arrangements \cite{OT}
states that the number of bounded regions of a real arrangement
equals the Euler characteristic of the complement of
its complexification. The Euler characteristic is 
$\,(-1)^{d+1} \mu(0)$ where $\mu(0)$ is the M\"obius number of
the intersection lattice, and, by Theorem \ref{thm-linear}, this is precisely the ML degree of
the associated linear model. Here we extend this relationship
between topology and the ML degree 
to statistical models which are given by nonlinear polynomials~$f_i$.
Working in the general setting of Section 2,
we shall prove the following:

\begin{thm} \label{euler}
Let $X$ be a smooth complete algebraic variety over $\cc$ of dimension
$d$, and let $D$
be the reduced divisor associated to $f = f_1^{u_1} \cdots f_n^{u_n}$.
Assume that the hypotheses (a), (b) and (c) below hold.
Then the ML degree equals
   $(-1)^d e_{top} (X \backslash D)$ where $e_{top}$
is the topological Euler
characteristic.
\end{thm}

Invoking Hironaka's theorem on resolution of singularities,
we fix a blow up $\pi : \tilde{X} \ra X$ such that
$\tilde{X}$ is smooth, and the
rational function $f$ pulls back to a proper morphism $\tilde{f}$ 
onto $\PP^1_{\C}$. The three {\bf hypotheses} are as follows:
 
\begin{enumerate}
\item[(a)] The inverse image $D' : = \pi^{-1}(D)$ of the divisor 
$D$ can be written as
$ D' = \tilde{D} + D_H$ where $ \tilde{D}$ is the support of the divisor 
${\rm div}(\tilde{f})$ of the pullback
$\tilde{f}$ of the rational function $f$, while $D_H$ is the 
horizontal divisor
consisting of the sum of all components of $D'$ which map onto $\PP_\C^1 $.
\item[(b)] The restriction of $\tilde{f}$ to $D_H \backslash \tilde{D}$ 
is a topological fiber bundle
over $\, \C^* \, = \,  \PP^1_{\C} \backslash \{0,  \infty \}$.
\item[(c)] The number of critical points of $\tilde{f}$ 
 on  $\tilde{X} \backslash \tilde{D}$ is finite.
\end{enumerate}


\begin{rem} \label{wowo}
Hypothesis (a) is crucial and depends on the exponents 
$u_i$. For instance, consider the
rational function $f = (y-x^2)(y+x^2)^{-1}$ on
$X = \PP^2_{\C}$. We get $\tilde{X}$ by blowing
up the origin twice. The  exceptional curve of the first blow-up 
belongs to the fiber $\{\tilde{f}=1 \}$ and hence is not
supported on ${\rm div}(\tilde{f})$. Hypothesis 
(a) is not satisfied for this example. If we take instead   
$f = (y-x^2)(y+x^2)^{-2} $ then hypothesis (a) is satisfied because the
exceptional curve belongs to the fiber $\{ f= \infty \}$
and is hence supported on ${\rm div}(\tilde{f})$,

Hypothesis (b) implies that the cohomology 
ranks of $D_H \backslash \tilde{D}$ can
be computed from that of $\cc^*$ and the fibers using
K\"unneth's formula. 
The alternating sum of the ranks is zero for the fibers,
and we get $\,e_{top}(D_H \backslash \tilde{D}) = 0$. In fact, the 
hypothesis (b) could be replaced by the more general condition
$\,e_{top}(D_H \backslash \tilde{D}) = 0$.  

Any proper map $f : X \to \mathbb{C}^{\ast}$ for $X$ smooth is a topological
fiber bundle if it is a submersion, i.e.  $df \neq 0$ for all points
in $X$. Therefore to check hypothesis (b), we need only find a controlled
stratification (in the sense of Thom-Mather theory \cite{mather})
of $D_H$ into locally closed smooth sets such that for
all points on each strata $S$, $d(\tilde{f}|_S) \neq 0$.  In
particular, this last condition will imply that the critical points of
$\tilde{f}$ on $\tilde{X} \backslash \tilde{D}$ are the same as the
critical points of $f$ on $X \backslash D$.
\end{rem}

\smallskip {\em Proof of Theorem \ref{euler}:}
Our method follows closely the proof of the Lefschetz hyperplane
section theorem (cf.~\cite{a-f1, a-f2}).
Moreover, since the complement is not necessarily compact we shall
use Borel-Moore homology \cite{BM} (see also \cite{bredon}).
We note that for compact spaces $X$ the  ordinary homology groups coincide
with the Borel-Moore homology groups. In the Borel-Moore
homology theory we have the following useful exact sequence to
be used below: if $X$ is locally compact,
$F$ is closed in $X$, and $ \,U : = X \backslash F$, we have
\begin{equation}
\label{bmsequence}
  \dots  \ra H_i(F) \ra H_i(X) \ra H_i (U) \ra \dots .
  \end{equation}
 Thus in this situation the Borel-Moore Euler characteristic is additive:
  \begin{equation}
\label{additive}
  e_{BM}(X) \quad = \quad e_{BM}(F) \, +\, e_{BM}(U). 
   \end{equation} 
  Finally, if $U$ is an even-dimensional orientable  manifold
  then Poincar\'e duality holds between
Borel Moore homology and ordinary cohomology, and 
$e_{BM}(U)$ coincides with the topological Euler number $e_{top} (U)$.
In our situation we get
$$
e_{top} (X \backslash D) \,\,  =  \,\, e_{top} (\tilde{X} 
\backslash D') \,\, = \,\,
  e_{BM} (\tilde{X} \backslash D')
 \,\,\, = \,\,\, e_{BM} (\tilde{X} \backslash \tilde{D})
\,-\, e_{BM} (D' \backslash \tilde{D} ). $$
The last equation follows from (\ref{additive}).
Hypothesis (b) implies $\, e_{BM} (D' \backslash \tilde{D} )  = 0\,$ 
 (see Remark \ref{wowo}), and hence it suffices
 to show that the ML degree equals
 $$  e_{top} (X \backslash D) \,\,  =  \,\,
 e_{BM} (\tilde{X} \backslash \tilde{D}) \,\,\, = \,\,\,
e_{top} (\tilde{X} \backslash \tilde{D}).$$
In other words, we may now simply erase the tilde 
and consider the case when $X$ is smooth and $f$ defines a proper morphism
 $X \backslash D \rightarrow \C^*$.

Let $\mathcal{C}$ denote the set of critical points of $f$ on $X$.
By hypothesis (c), this set is finite and the ML degree equals
its cardinality counting multiplicities:
 $$ \mu \quad = \quad \sum_{p \in \mathcal{C}} \mu_p . $$
  The multiplicity $\mu_p$ of a critical point $p$ of $f$ is known as
 the {\em Milnor number} at $p$ of the
hypersurface $F_p = \{ x \in X \,: \, 
f(x) = f(p)\}$. Milnor  \cite{Mil2} showed that this algebraic invariant
of a singularity
has the following topological interpretation. Consider a coordinate
chart around the point $p$ and intersect the fiber
$ \,F_{\epsilon} := \{ x  \, |\, f(x) = f(p) + \epsilon\}\,$
with a  ball of radius $\delta$ around $p$.
For $ \epsilon \ll \delta$ this intersection is the
{\em Milnor fiber}. Milnor \cite{Mil2} showed that the Milnor fiber is
homotopy equivalent to a bouquet of $\mu_p$ spheres of dimension $d-1$.

Each singular fiber  is obtained  (cf.~\cite{BPV, Mil2}) from a smooth
fiber  by replacing the Milnor
fiber  by a contractible set. The Borel-Moore exact sequence 
(\ref{bmsequence}) implies that
the Euler number of a singular fiber $F'$ is obtained from the
Euler number of a smooth fiber $F$ by adding $ - (-1)^{d-1} \Sigma_{p 
\in F'} \mu_p$.

Then the Euler number of the union of the singular fibers equals
$|\mathcal{C}|$ times the Euler number of a smooth fiber $F$ plus the
correction $ - (-1)^{d-1} \mu$. Applying K\"unneth's formula to the 
fiber bundle defined by $f$ on $X \backslash D$ minus the union of
the singular fibers, and then applying the additivity formula
(\ref{additive}), we conclude that $\,e_{BM}(X \backslash D) \, =\,
e_{top}(X \backslash D) \,= \,(-1)^d \mu \,$ as desired.  $\Box$

\smallskip

\begin{ex}
 For another illustration  consider Example  \ref{4 curves} with $X = \PP^2_{\C}$.
The generic ML degree was $9$ but it decreased  by $4$
when $ u =  (u_1,u_2,u_3,u_4)$ is  a general solution of
 $u_1 + u_2 + 2 u_3 + 2 u_4 = 0$. This is consistent with
Theorem \ref{euler}, 
because for such $u$ the divisor $D$ loses the component at infinity.
The difference is a  projective line minus $6$ points, which has Euler number $-4$.
Consider our hypotheses when  $\tilde X$ is the blow-up 
of $X$ at the origin. If $ u_1 + u_2 +  u_3 +  u_4 \neq 0 $ then
the exceptional curve is part of $\tilde{D}$ and
Theorem  \ref{euler} is valid. On other hand, if
$ u_1 + u_2 +  u_3 +  u_4 = 0 $ then it maps to $\PP^1$
under a rational map of degree $\geq 2$,
so hypothesis (a) does not hold.
The philosophy of this example 
is that, even if the divisor $D$ is locally biholomorphic to an
arrangement of  hyperplanes,
genericity of the exponents $u_i$ may be necessary for the 
topological formula
of Theorem \ref{euler} to hold. 
\qed
\end{ex}

Theorems \ref{thm1},  \ref{toricformula3} and \ref{toricformula4}
 offer combinatorial formulae for the ML degree and
 hence (using Theorem \ref{euler}) for the 
  Euler number of the complement
 of an arrangement of generic hypersurfaces
 $\{ f_i  = 0\}$  in $X = \PP^d_{\C}$ or in a toric manifold.
In each of these theorems, the
 combinatorial number becomes an upper bound
 for the ML degree when the coefficients of
 the $f_i$'s are special.   This semi-continuity
 principle will be explained by the
 following general topological result,
 in which also the  underlying manifold is
 allowed to vary.
 
 \begin{thm} \label{semiconti} Assume we are given a 
 one-parameter smooth proper family $X_t$ of
complex manifolds
over the unit disk $B : = \{ \,t \in \C \,:\, |t| < 1 \}$, and a family of rational
functions $f_t$ on $X_t$, such that
\begin{enumerate}
\item for $\, t \neq 0 \,$ the divisor $D(t) $ defined by $f_t$ has GNC, and
\item for $\, t=0 \,$ the divisor $D(0)$ defined by $f_0$ has the same 
homology class as $D(t)$ for the natural
differentiable trivialization of the family $X_t$.
\end{enumerate}
Then the ML degree of $f_0$ is less than or equal to the ML degree of
$f_t$.
\end{thm}

In order to understand the second hypothesis, let us recall 
Ehresmann's Theorem \cite{Ehre}: any proper submersion $\phi : \X \ra B$ of 
differentiable
manifolds is a differentiable fiber bundle, i.e., if $U$
is a sufficiently small open set in  $B$, there is a local 
diffeomorphism between  $\phi^{-1}(U)$ and $U \times F$,
for  a fiber $F$, and this diffeomorphism is 
compatible with the two
projections to $U$.

\medskip \noindent {\sl Sketch of Proof.}
Let $\X \ra B $ be the total space of the family $\{X_t\}_{t \in B}$.
We consider the function
$\phi : \X \ra B \times \PP^1$, given by $\phi (x) = (t(x), f(x))$.

Consider the locus $\Xi \subset \X$ given by the vanishing of the
vertical differential of $\phi$: this is  the local complete intersection
defined by the $d$ partials $\partial f / \partial x_i = 0$,
where $x_1, \ldots, x_d$ are local coordinates on the fibers
provided by the Implicit Function Theorem.
At each critical point $p$ of $f_0$,
the locus $\Xi$ has dimension $1$, and thus, in a
neighborhood of $p$, the morphism $\Xi \ra B$ is finite, whence its
degree is locally constant.
This establishes the desired semi-continuity.
\qed

\medskip

In general, it might be difficult to show that a given rational function
$f_0$ has a perturbation as above, or we might want to calculate
the ML degree with more algebraic precision. The results in the
next section may help.

\section{Logarithmic vector fields}

In this section we will show that the  formula for the ML degree given
by the logarithmic Chern number (Theorem \ref{chernclassthm})
holds in greater generality. Returning to the setting of
 Chapter 2, we consider logarithmic vector fields along $D$.
  Again, our definition differs slightly from the one given by Saito \cite{Sai}.

\begin{df} If $D$ is a  reduced divisor on a factorial variety $X$,
the {\em sheaf of logarithmic vector fields} $\Theta_X(- log D)$ 
is the dual $ \HHH om_{\hol_X} (\Omega^1_X (log D), \hol_X )$.
\end{df}

Recall that the {\em tangent sheaf} $\Theta_X$ is 
$ \HHH om_{\hol_X} (\Omega^1_X , \hol_X )$, the dual of the 
$1$-forms on $X$. The inclusion of 
$\Omega^1_X$ into $\Omega^1_X (log D)$,
studied in Lemma~\ref{exactlem},
dualizes to an inclusion of the logarithmic
vector fields into the tangent sheaf.

\begin{prop} \label{moresheaves}
We have the following exact sequence of sheaves on $X$:
$$ 0 \ra \Theta_X(- log D) \ra \Theta_X \ra \bigoplus_{i=1}^r \hol_{D_i}
(D_i) \ra  \EE xt^1_{\hol_X} (\Omega^1_X (log D), \hol_X ). $$
If $X$ is smooth, then the  rightmost homomorphism is onto,
and  the total Chern class of the sheaf of logarithmic
vector fields equals
\begin{equation}
\label{ChernOfTheta}
c_{tot} \bigl(\Theta_X(- log D)\bigr) \quad = \quad
\frac{ c_{tot} (\Theta_X)  \cdot c_{tot} (\EE xt^1_{\hol_X} (\Omega^1_X (log D), \hol_X ))}
{\Pi_{i=1}^{r} ( 1 + z D_i)}.
\end{equation}
\end{prop}

\begin{proof}
Dualizing the sequence
$0 \ra \hol_X(-D_i) \ra \hol_{X} \ra \hol_{D_i} \ra 0$,
we get
$$ \HHH om_{\hol_X} (\hol_{D_i}, \hol_X ) \,\,\, = \,\,\, 0 \quad
\hbox{and} \quad 
  \EE xt^1_{\hol_X} (\hol_{D_i}, \hol_X ) \,\,\, \cong \,\,\,
  \hol_{D_i} (D_i) . $$
 Hence the $\EE xt_{\hol_X}$-sequence
 gotten by dualizing the sequence (\ref{lemmasequence})
 has the form
 \begin{eqnarray*}
& 0 \,\,\, \ra \,\, \,\,  0 \, \ra \, \Theta_X(- log D) \, \ra \, \Theta_X \ra
 \qquad \qquad \qquad \qquad \qquad \\ & \quad
 \bigoplus_{i=1}^r \hol_{D_i} (D_i) 
 \ra  \EE xt^1_{\hol_X} (\Omega^1_X (log D), \hol_X ) \ra  
 \EE xt^1_{\hol_X} (\Omega^1_X , \hol_X ) \ra \cdots
 \end{eqnarray*}
 This is the first statement of Proposition \ref{moresheaves}.
 If $X$ is smooth then the cotangent sheaf $\Omega^1_X$ is free, and
 we have  $\EE xt^1_{\hol_X} (\Omega^1_X , \hol_X ) =0$.
 The formula (\ref{ChernOfTheta}) follows from the  
multiplicativity of the total Chern class.
\end{proof}

\begin{rem} \label{logvecchara}
The two leftmost maps in the exact sequence of  Proposition  \ref{moresheaves}
characterize the logarithmic vector fields on $X$ along $D$ as those vector fields
$\xi \in \Theta_X$ which  satisfy $\,\xi (F_i) \equiv 0
\, ({\rm mod} \, F_i)\,$ for all $i$.
In other words, for each $ i=1, \ldots, h$, the vector field
$\xi = \Sigma_{j=1}^d  \xi_j \partial / \partial x_j$ has the property
that there exist functions $\psi_i$ such that
$\xi (F_i) : = \Sigma_{j=1}^d  \xi_j \partial  F_i / \partial x_j
= \psi_i F_i$.
\end{rem}

\begin{rem}
An interesting case will be the one where $X$ is smooth and
the sheaf $\EE xt^1_{\hol_X} (\Omega^1_X (log D), \hol_X )$ has
zero-dimensional support and length $p$. Here the top Chern
class of $\Theta_X(- log D)$ and the top Chern class
of $\Omega^1_X (log D)$ differ by $(-1)^{d-1} p$,
and the count in part 3 of Theorem \ref{chernclassthm}
changes accordingly.
\end{rem}

As in the proof of Theorem \ref{chernclassthm},
set $\,\sigma :=  dlog (f)\,$ and let 
$Z_\sigma $ be the subscheme
of $X$ defined by the vanishing of $\sigma$.
The restriction of $Z_\sigma$ to the open set
$V = X \backslash D$ is the critical locus of $f$.
The reason why the sheaf $\Theta_X(- log D)$ is important is that
it enters directly into the algebraic description
of $Z_\sigma$.

Indeed, the section $\sigma$ of 
$\,\Omega^1_X (log D)\,$
corresponds to an exact sequence
$$ 0 \ra \hol_X \ra \Omega^1_X (log D) \ra \EE \ra 0 ,$$
where $ \hol_{Z_{\sigma}}$ is the kernel of
$ \EE xt^1_{\hol_X} (\,\EE , \hol_X ) \ra
\EE xt^1_{\hol_X} (\,\Omega^1_X (log D), \hol_X )$.
Hence the $\EE xt_{\hol_X}$-sequence
 gotten by dualizing the previous sequence  has the form
\begin{equation}
\label{Zsequence}
 0 \ra   \HHH om_{\hol_X} (\, \EE , \hol_X )
\ra  \Theta_X(- log D) \ra \hol_X  \ra
\hol_{Z_{\sigma}} \ra 0.
\end{equation}
This leads to the following more refined
formula for the ML degree.

\begin{thm} \label{thetatheorem}
Assume that  $X$ is smooth, $D$ meets every curve,
the sheaf $\Theta_X(- log D)$ is locally free, and 
$Z_{\sigma}$ does not intersect the  divisor $D$.
Then the number of critical points of $f$ on $ V = X \backslash D$ equals
$ (-1)^d c_d ( \Theta_X(- log D) )$.
\end{thm}
\begin{proof}
Using (\ref{Zsequence}),
we can view $Z_{\sigma}$ as the locus of zeros of the locally free
sheaf dual to $\Theta_X(- log D)$, whose top Chern class is
$ (-1)^d c_d ( \Theta_X(- log D) )$.
\end{proof}

\begin{rem}
The sequence (\ref{Zsequence}) gives the following
description of the equations defining the critical points.
The ideal sheaf of the subscheme $Z_\sigma$ is generated
by the functions $ \Sigma_{i=1}^r u_i \psi_i$, where $\xi$ varies
among the logarithmic vector fields and $(\psi_1,\ldots,\psi_r)$
is derived from $\xi$ as in Remark \ref{logvecchara}.
This holds because  the  section $\sigma = dlog(f) $
factors through the homomorphism $\hol_X^r \ra \Omega^1_X ( log D)$,
thus the homomorphism dual to the section $\sigma$ also factors through the dual
homomorphism $\,\Theta_X(- log D) \ra \hol_X^r,
\, \xi \mapsto (\psi_1,\ldots,\psi_r) $.
\end{rem}

We next present two examples which illustrate the hypotheses
of Theorem \ref{thetatheorem}.  Here $X$ is a smooth surface
(i.e.~$ d = 2$) with local coordinates $x$ and $y$.

\begin{ex} Let $h = 2$ with  $F_1 = x$ and $F_2 =
x^v - y^u$. A vector field $\xi = a (x,y) \cdot \partial  / \partial x +
b(x,y) \cdot \partial  / \partial y$ is logarithmic if and only if
 $ a =  x \psi_1 $ and $ b = (v/u) \ y \ \psi_1 + \lambda (x^v - y^u)$,
for some function $\lambda $. Observe that $\psi_2 = v \psi_1 - y ^{u-1} \ \lambda$.
We conclude that the sheaf $\Theta_X(- log D)$ is locally free of rank~$2$.

The origin does not belong to the subscheme $Z_{\sigma}$ if there is a function
$u_1 \psi_1 + u_2 \psi_2 = (u_1 + v \cdot u_2 ) \psi_1 - u_2 y ^{u-1}
\ \lambda$, for some $ \psi_1 $ and $\lambda$,
which does not vanish at the origin.
This holds if and only if $\, u_1 + v \cdot u_2 \,\neq \, 0$. \qed
\end{ex}

\begin{ex}
Let $h = 3$ with $F_1 = x, F_2 =y$ and $F_3 = x-y$. 
Logarithmic vector fields $\xi$ have the form
$\, x \psi_1 \cdot \partial  / \partial x  + y \psi_2 \cdot \partial  / \partial y \,$
with $\,  x \psi_1 - y \psi_2 \,$ divisible by $x-y$. This implies
 $ \psi_1 = \lambda  + y \cdot  \mu$,  $ \psi_2 = \lambda  + x \cdot  \mu$,
and $  \psi_3 =  \lambda$, for any functions
$\lambda , \mu$. Thus $\Theta_X(- log D)$ is locally free of rank two.
The origin does not belong to the subscheme $Z_{\sigma}$ if and only if
$\, u_1 + u_2 + u_3 \neq 0$.
\end{ex}

\medskip

\noindent {\bf Acknowledgements:}
This project grew out of discussions with
Lior Pachter at the December 2003 workshop
on {\em Computational Algebraic Statistics}
at the American Institute for Mathematics (AIM)
at Palo Alto. Many thanks to both Lior and AIM for getting us started,
and many thanks to MSRI Berkeley and the organizers
of the Spring 2004 MSRI program
{\em Topology of Real Algebraic Varieties},
which gave the four of us the opportunity to collaborate.
We also thank Oleg Viro for letting us include Theorem \ref{thm:viro}.
Amit Khetan was supported by an NSF postdoctoral fellowship (DMS-0303292).
Bernd Sturmfels was supported by
the Hewlett Packard Visiting Research Professorship 2003/2004
at MSRI Berkeley and in part  by the NSF (DMS-0200729).

\bigskip

\noindent {\bf Authors' addresses:}

\medskip

\noindent  Fabrizio Catanese, Lehrstuhl Mathematik VIII,
Universit\"at Bayreuth, NWII,  D-95440 Bayreuth, Germany,
{\tt Fabrizio.Catanese@uni-bayreuth.de}
\medskip

\noindent  Serkan Ho\c sten, Department of Mathematics,
San Francisco State University, San Francisco, CA 94132, USA,
{\tt serkan@math.sfsu.edu}

\medskip

\noindent Amit Khetan, Department of Mathematics,
University of Massachusetts, Amherst, MA 01002, USA,
{\tt khetan@math.umass.edu}

\medskip

\noindent Bernd Sturmfels, Department of Mathematics,
University of California, \break Berkeley, CA 94720, USA,
{\tt bernd@math.berkeley.edu}

\end{document}